%% file: fb1.tex
\documentclass[a4paper,12pt]{article}
\usepackage[utf8]{inputenc}
\usepackage{amsmath}
\usepackage{amssymb}
\usepackage{euler}
\usepackage{eucal}

\usepackage[T1]{fontenc}
\usepackage[english]{babel}

\usepackage{tikz} \usetikzlibrary{positioning,shapes}
\usepackage{pgfplots} \usetikzlibrary{plotmarks}

\usepackage[unicode=true,plainpages=false]{hyperref}
\usepackage{graphicx}
\hypersetup{colorlinks=false,pdfborder={0 0 0.1},linkcolor=magenta,anchorcolor=magenta,urlcolor=blue,citecolor=blue,
          }

\usepackage{algorithmic}
\usepackage[ruled]{algorithm}
\usepackage{amsthm}
\theoremstyle{definition}

\newtheorem{theorem}{Theorem}

\let\eps=\varepsilon

\let\le=\leqslant
\let\ge=\geqslant

\def\O{\mathcal{O}}

\title{Low-rank approximation in the numerical modeling of the Farley-Buneman instability in ionospheric plasma \thanks{This work was partially supported by RFBR grants
13-01-12061 (ofi-m-2013), 12-01-91333 (nnio-a), 11-01-00549-a, 12-01-33013 (mol-a-ved), 12-01-31056, Government Contracts 16.740.12.0727, $\Pi$1112, 8500
and the support programmes of the RAS Presidium and RAS Department of Mathematical Sciences.}}

\author{S. V. Dolgov\thanks{Max-Planck-Institut f{\"u}r Mathematik in den Naturwissenschaften, Inselstr. 22-26, D-04103 Leipzig, Germany. Institute of Numerical Mathematics of Russian Academy of Sciences, Moscow, Russia ({\tt sergey.v.dolgov@gmail.com})},
A. P. Smirnov\thanks{Lomonosov Moscow State University, Moscow, Russia ({\tt sap@cs.msu.su})},
E. E. Tyrtyshnikov\thanks{Institute of Numerical Mathematics of Russian Academy of Sciences, Moscow, Russia.
Lomonosov Moscow State University, Moscow, Russia.
Moscow Institute of Physics and Technology, Moscow, Russia.
University of Podlasie, Siedlce, Poland (\tt eugene.tyrtyshnikov@gmail.com)}}

\date{August 27, 2013}

\begin{document}
\bibliographystyle{siam}

\maketitle

\begin{abstract}
We consider the numerical modeling of the Farley-Buneman instability development in the earth's ionosphere plasma.
The ion behavior is governed by the kinetic Landau equation in the four-dimensional phase space, and since the finite difference discretization on a tensor product grid is used, this equation becomes the most computationally challenging part of the scheme.
To relax the complexity and memory consumption, an adaptive model reduction using the low-rank separation of variables, namely the Tensor Train format, is employed.

The approach was verified via the prototype MATLAB implementation.
Numerical experiments demonstrate the possibility of efficient separation of space and velocity variables, resulting in the solution storage reduction by a factor of order tens.



\end{abstract}

\section{Introduction}

The Farley-Buneman (FB) instability arises in weakly ionized plasma of E-region of Earth's ionosphere.
This instability is generated in the plasma with magnetized electrons and unmagnetized ions in the electric field  which is directed perpendicular to the geomagnetic field  \cite{dimant-fb-linear-2004}.
In the E-region electrons are magnetized by the geomagnetic field, while ions are not magnetized due to frequent collisions with neutral gas particles.
As a result, the velocity distribution of electrons is shifted relative to that of ions by an electrical drift velocity.
Appropriate conditions for the onset of instability occur in the equatorial and polar zones of the E-region of the Earth's ionosphere, where the instability manifests itself as low-frequency plasma oscillations with wavelengths in the scale of meters.

The first papers which investigated the FB instability were published by Farley \cite{farley-plasma-1963} and Buneman \cite{buneman-plasma-1963} independently.
They used a linear theory.
Using the kinetic equations, Farley showed that the strong external electric field leads to the instability and appearance of waves in plasma.
Buneman obtained the dispersion relation from the fluid theory, which indicates that the growth of the instability is possible only when the electron drift velocity exceeds some threshold.
It means that the external electric field should be large enough.

The linear theory permits to derive the threshold conditions, giving the necessary conditions for instability development, but this theory can not describe the process of instability saturation.
The latter can be described only on the base of a nonlinear theory, which was developed for a long time \cite{skadron-fb-nonlin-1969,sudan-electrojet-1973,hamza-fb-1995}, but still has a limited application.

The nonlinear models of the FB instability are based on several nonlinear two- and three-dimensional partial differential equations.
The quantitative solution of these FB instability equations needs to conduct computer modeling, by which one can evaluate significance and applicability of those or other theories.

First computer simulations of FB instabilities were based on the fluid theory \cite{newman-nonlinear-jet-1981}.
With the development of computers, more advanced models and methods have been proposed.
The particle method was used in \cite{machida-fb-1988,schlegel-pic-fb-1994,janhunen-fb-1994,dimant-fb-pic-2004}, and the combined method based on particles and fluid equations was exploited in \cite{oppenheim-fb-sim-1996,oppenheim-fb-sat-1996,dyrud-fb-rocket-2006}.
The usage of a fluid model for electrons and ions leads to a nonphysical result: the growth rate of the instability increases infinitely with the wave number.

Contrarily, the kinetic Landau damping allows to describe the process of the FB instability saturation.
While both electrons and ions are prone to Landau damping resulting in a suppression of the instability, electron Landau damping is only effective at a short-wavelength, high-frequency range,
but the wave growth rate is restrained efficiently by the ion Landau damping only.
This allows to employ the fluid model for electrons.

In this article we will use a hybrid model for the FB instability proposed in \cite{kovalev-plasma-2008,kovalev-plasma2Te-2009,kovalev-plasma3kinetic-2009,kovalev-plasma4effect-2009}.
The model is based on the following equations:  2D fluid equation for the electron density, 4D kinetic equation for ions and 2D Poison equation for the electric field potential.
The growth and saturation of the FB process in the hybrid model was studied numerically \cite{kovalev-plasma-2008}, using the derivatives discretization on multi-dimensional meshes.
Most of the computer time was spent for the numerical solution of the ion kinetic equation on a 4D mesh in the phase space (2D-space, 2D-velocity).
The size of the 4D array representing the kinetic solution is about $10^9-10^{12}$ bytes, and such a large amount of RAM memory, as well as the corresponding number of computer operations in each time step insist to use high-performance parallel systems.
The FB instability simulation described in \cite{kovalev-plasma-2008} was implemented in such a parallel code for the supercomputer Blue-Jean P.

In this work a new tensor approximation approach is used for the solution of the algebraic system obtained from the numerical approximation of the 4D kinetic equation.
The main idea is in separation of spatial and velocity coordinates of the phase space, and approximation of the whole function by a proper low-term sum of direct products of spatial and velocity contributions (see Section \ref{sec:tt} for details).
Equipped with efficient variational techniques for adaptive approximations directly in the separated form, the tensor structuring provided the reduction of the necessary RAM memory in 20 times, and the computational time at least twice.
This simulation was conducted via a sequential prototype MATLAB code at PC.
More significant performance gain may be expected after refactoring both the program code into the parallel optimized version, and the discretization schemes to high-order accurate ones, which is a matter of a forthcoming work.

The paper is organized as follows.
In Section \ref{sec:fb-setting} we formulate the hybrid model of the FB instability, and discuss the output quantities of interest.
The main Section \ref{sec:solver} is devoted to the spatial and time discretization schemes, tensor product approximation formats and methods, and how to cast the solution scheme steps into the tensor structures.
In Section \ref{sec:num}, the numerical simulation is presented, which compares the initial hybrid model where the ion distributions are stored straightforwardly as 4D arrays with the tensor approximation counterpart.
Finally, Section \ref{sec:conclusion} gives concluding remarks and points out perspectives.

\section{Problem statement}\label{sec:fb-setting}
The Farley-Buneman instability occurs in the so-called E-region of ionosphere (90-100 km),
where the electrons are affected by the geomagnetic field, but the ions are not \cite{farley-plasma-1963,buneman-plasma-1963}.
Though a very general kinetic behavior of a plasma is accurately described by the Vlasov-Fokker-Planck equation,
the E-region conditions allow to consider a simplified hydrodynamical model of electrons.

The coordinate axes are introduced along the electric and magnetic fields, so that
$$
\begin{array}{rcl}
\mathbf{B}_0 & = & \begin{bmatrix}0 & 0 & B_0 \end{bmatrix} \\
\mathbf{E}_0 & = & \begin{bmatrix}0 & E_0 & 0 \end{bmatrix} \\
\mathbf{V}_0 & = & \begin{bmatrix}V_0 & 0 & 0 \end{bmatrix} = \mathbf{E}_0 \times \frac{\mathbf{B}_0}{B_0^2}, \\
\end{array}
$$
where $\mathbf{B}_0$ is the geomagnetic field, $\mathbf{E}_0$ is an external electric field, and $\mathbf{V}_0$ is a drift velocity.
We use the two-dimensional model in the $x,y$ plane only, since the time scale of the processes along $\mathbf{B}_0$ is significantly smaller, as well as their magnitude.

\subsection{Electron model}
The behavior of electrons is assumed to be governed by the standard continuity equations:
\begin{equation}
\begin{array}{rcl}
 \dfrac{\partial n_e}{\partial t} & = & - \nabla (n_e \mathbf{V}_e), \\
 m_e \dfrac{d \mathbf{V}_e}{dt} & = & -e (\mathbf{E}_0 - \nabla \Phi + \mathbf{V}_e \times \mathbf{B}_0) - \dfrac{\nabla(n_e T_e)}{n_e} - m_e \mathbf{V}_e \nu_{en},
\end{array}
\label{eq:el_cont}
\end{equation}
where $n_e = n_e(x,y)$ is the electron concentration, $\mathbf{V}_e$ is the total electron velocity, $\Phi$ is an electric potential,
$T_e$ is an electron temperature (the effect of electron heating was studied in \cite{kovalev-plasma2Te-2009}), $\nu_{en}$ is the average electron-neutral collision rate, and $m_e$ and $e$ are the electron mass and charge, respectively.

We will normalize all quantities to the characteristic scales of the model as shown in Table \ref{tab:rescale}.
\begin{table}[h!]
\centering
\caption{Typical model scales and renormalizations of quantities}
\label{tab:rescale}
\begin{tabular}{|cc|cc|c|} \hline
\multicolumn{2}{|c|}{Initial Quantity} & \multicolumn{2}{|c|}{Scale} & Dimensionless quantity \\ \hline
Time & $t$ & $1/\nu_{in}$ & [s] & $t:=t \nu_{in}$ \\ \hline
Space(x) & $x$ & $l = v_{T_i}/\nu_{in}$ & [m] & $x:=x/l$ \\ \hline
Space(y) & $y$ & $l = v_{T_i}/\nu_{in}$ & [m] & $y:=y/l$ \\ \hline
Velocity(x) & $v$ & $v_{T_i}=\sqrt{T_i/m_i}$ & [m/s] & $v:=v/v_{T_i}$ \\ \hline
Velocity(y) & $w$ & $v_{T_i}=\sqrt{T_i/m_i}$ & [m/s] & $w:=w/v_{T_i}$ \\ \hline
Temperature(e) & $T_e$ & $T_i$ & [J] & $T_e:=T_e/T_i$ \\ \hline
El. potential & $\Phi$ & $T_i/e$ & [V] & $\phi = e \Phi/T_i$ \\ \hline
\end{tabular}
\end{table}
In addition, introduce the following aggregated dimensionless quantities:
\begin{equation}
\psi = \dfrac{\nu_{en} \nu_{in}}{\Omega_e \Omega_i}, \qquad \theta = \left(\dfrac{m_e \nu_{en}}{m_i \nu_{in}}\right)^{1/2},
\end{equation}
where $\nu_{in}$ is the average ion-neutral collision rate, $\Omega_{e,i}$ are the cyclotron frequencies of the electrons and ions, respectively, and $m_i$ is the ion mass.
Since in the E-region it holds $\omega \ll \nu_{en}$, where $\omega$ is the plasma frequency, we can neglect the electron inertia in \eqref{eq:el_cont}.
Taking all the considerations presented above together, the electron equation can be written as follows,
\begin{equation}
\begin{array}{rcl}
 \dfrac{1}{\psi\sqrt{T_e}} \dfrac{\partial n_e}{\partial t} & = & \Delta (T_e n_e) + \dfrac{\partial}{\partial x} \left(\dfrac{V_0}{v_{T_i}\psi\sqrt{T_e}} - \dfrac{1}{\theta\sqrt{T_e\psi}} \dfrac{\partial \phi}{\partial y} - \dfrac{\partial \phi}{\partial x} \right) n_e \\
 & + & \dfrac{\partial}{\partial y} \left(\dfrac{eE_0v_{T_i}}{T_i \nu_{in}} + \dfrac{1}{\theta\sqrt{T_e\psi}} \dfrac{\partial \phi}{\partial x} - \dfrac{\partial \phi}{\partial y} \right) n_e,
\end{array}
\label{eq:electron}
\end{equation}
where $\phi$ is a dimensionless electric potential due to the non-uniform distribution of electron and ion densities (see Table \ref{tab:rescale}), and $T_i$ is the initial temperature of ions.
Note that the temperatures are measured in Joules, i.e. $T_* = k \hat T_*$  if $\hat T_*$ is measured in Kelvins.

Since on a large scale we assume the spatial uniformity, Equation \eqref{eq:electron} is posed in a box $[0,L]^2$ with the periodic boundary conditions.

\subsection{Poisson equation for the electric potential}
Given the charge distribution $\rho$, the electric potential obeys the Poisson equation,
$$
\Delta \Phi = \dfrac{1}{\varepsilon_0} \rho.
$$
In our case, $\rho$ appears due to the non-uniformity of the particle densities, $\rho = e(n_e-n_i)$.
Therefore, after the renormalization, obtain the final equation
\begin{equation}
\Delta \phi = \dfrac{e^2}{\varepsilon_0 m_i \nu_{in}^2}(n_e-n_i).
\label{eq:poisson}
\end{equation}
As previously, we pose periodic boundary conditions on $[0,L]^2$.

\subsection{Kinetic description of ions}
Since the electron model is already simplified to a hydrodynamical description, we do not need
to employ the true Landau collision operator.
Instead, we use the so-called Bhatnagar-Gross-Krook (BGK) relaxation term for the velocity distribution of ions \cite{kovalev-plasma3kinetic-2009}.
In such a way simplified Vlasov-Fokker-Planck equation reads
$$
\dfrac{\partial f(x,y,v,w,t)}{\partial t} + \mathbf{v} \cdot \nabla_{x,y} f + \dfrac{e(\mathbf{E}_0-\nabla \Phi)}{m_i} \cdot \nabla_{v,w} f = -\nu_{in} (f-f_0),
$$
where $f_0$ is the distribution of the neutral particles, which is assumed to be Maxwellian, and $\mathbf{v}=\begin{bmatrix}v & w\end{bmatrix}$ is the velocity vector corresponding to the space vector $\begin{bmatrix}x & y\end{bmatrix}$.
Finally, in the dimensionless quantities,
\begin{equation}
 \dfrac{\partial f}{\partial t} + v \dfrac{\partial f}{\partial x} + w \dfrac{\partial f}{\partial y} - \dfrac{\partial \phi}{\partial x} \dfrac{\partial f}{\partial v} + \left(\dfrac{e E_0}{m_i v_{T_i} \nu_{in}} - \dfrac{\partial \phi}{\partial y}\right) \dfrac{\partial f}{\partial w} = f_0-f,
\label{eq:vfp}
\end{equation}
where
\begin{equation}
f_0 = \dfrac{n_i}{2\pi} \exp\left(-\dfrac{v^2+w^2}{2}\right), \qquad n_i = \int\limits_{\mathbb{R}^2} f(x,y,v,w) dv dw.
\label{eq:f0}
\end{equation}

\subsection{Output}
There are several observable quantities that can be predicted by our model and verified experimentally.
First, the total electric field reads
$$
\mathbf{E}_{tot} = \mathbf{E}_0 + \nabla \Phi,
$$
and an \emph{additional} field $\mathbf{E}_{add} = \nabla \Phi$ appeared due to nonuniformly drifting charges can be considered as a perturbation to the initial \emph{signal} $\mathbf{E}_0$ (e.g. an ambient electric field).
Therefore, it is interesting to track the average magnitude of the additional field,
$$
E_{add} = \sqrt{(\partial \Phi/\partial x)^2 + (\partial \Phi/\partial y)^2} = \frac{T_i}{le} \sqrt{(\partial \phi/\partial x)^2 + (\partial \phi/\partial y)^2}
$$
(in the latter expression $x$ and $y$ are dimensionless).

Since the additional field varies in space and time, another characteristic is the principal wavelengths.
As was discussed in \cite{kovalev-plasma-2008,kovalev-plasma-2009}, despite the fact that the initial field $\mathbf{E}_0$ was directed along the $y$ axis, the drift vector rotates during the development of the Farley-Buneman process.
The spatial spectral intensities are computed as the squared values of the two-dimensional Fourier transformed field $\nabla \Phi$,
$$
\hat E^2(k_x,k_y) = |\hat E_x|^2 + |\hat E_y|^2, \quad \hat E_i(k_x,k_y) = \int E_i(x,y) \mathrm{e}^{-\mathrm{i} k_x x} \mathrm{e}^{-\mathrm{i} k_y y} dx dy.
$$

\section{Discretization and solution scheme}\label{sec:solver}
\subsection{Space and time discretizations}
For the spatial variables $x,y$, as well as for the velocities $v,w$ in \eqref{eq:vfp} we employ the
finite difference schemes, taking into account the periodic boundary conditions.
Though the velocity space is initially the whole plane, the actual distribution $f$ in \eqref{eq:vfp}, being a perturbed Maxwellian, decays rapidly with increasing $v$ and $w$, and we may shrink the domain to a cube $(v,w) \in [-v_{max},v_{max}]^2$, posing the periodic boundary conditions.

To treat the nonlinearity, we consider the following time splitting:
\begin{enumerate}
 \item Electron equation \eqref{eq:electron}, diffusion part.
 \item Electron equation \eqref{eq:electron}, advection part.
 \item Poisson equation \eqref{eq:poisson}.
 \item Ion equation \eqref{eq:vfp}, convections (sequentially in all variables).
 \item Ion equation \eqref{eq:vfp}, reaction (BGK relaxation).
\end{enumerate}
This is exactly the splitting scheme used in \cite{kovalev-plasma-2008,kovalev-plasma-2009}.
Unfortunately, it provides only the first order of accuracy.
However, it is not so difficult to develop a second-order accurate linearization.

Note from \eqref{eq:electron}, \eqref{eq:poisson} and \eqref{eq:vfp} that the global problem for $u = \begin{bmatrix}n_e & f\end{bmatrix}$ reads
$$
\dfrac{\partial u}{\partial t} = \left(\mathcal{A}_{ed}(u) + \mathcal{A}_{ec}(u) + \mathcal{A}_{ic}(u) + \mathcal{A}_{ir}(u)\right) u,
$$
where ``$ed$'' stands for ``electron diffusion'', ``$ec$'' for ``electron convection'', ``$ic$'' for ``ion convection'', and ``$ir$'' for ``ion reaction''.
Therefore, we may use established second-order splitting schemes for \emph{quasi-linear} problems, for example, the one proposed by Marchuk and Strang \cite{strang-split-1968,marchuk-split-1995}:
\begin{equation}
\begin{array}{rcl}
\tilde u & = & u^0 + \frac{\tau}{2} \left(\mathcal{A}_{ed}(u^0) + \mathcal{A}_{ec}(u^0) + \mathcal{A}_{ic}(u^0) + \mathcal{A}_{ir}(u^0)\right) u^0 \\
u^{1/8} & = & (I-\frac{\tau}{4} \mathcal{A}_{ic}(\tilde u))^{-1} (I+\frac{\tau}{4} \mathcal{A}_{ic}(\tilde u)) u^{0} \\
u^{2/8} & = & (I-\frac{\tau}{4} \mathcal{A}_{ir}(\tilde u))^{-1} (I+\frac{\tau}{4} \mathcal{A}_{ir}(\tilde u)) u^{1/8} \\
u^{3/8} & = & (I-\frac{\tau}{4} \mathcal{A}_{ed}(\tilde u))^{-1} (I+\frac{\tau}{4} \mathcal{A}_{ed}(\tilde u)) u^{2/8} \\
u^{4/8} & = & (I-\frac{\tau}{4} \mathcal{A}_{ec}(\tilde u))^{-1} (I+\frac{\tau}{4} \mathcal{A}_{ec}(\tilde u)) u^{3/8} \\
u^{5/8} & = & (I-\frac{\tau}{4} \mathcal{A}_{ec}(\tilde u))^{-1} (I+\frac{\tau}{4} \mathcal{A}_{ec}(\tilde u)) u^{4/8} \\
u^{6/8} & = & (I-\frac{\tau}{4} \mathcal{A}_{ed}(\tilde u))^{-1} (I+\frac{\tau}{4} \mathcal{A}_{ed}(\tilde u)) u^{5/8}, \\
u^{7/8} & = & (I-\frac{\tau}{4} \mathcal{A}_{ir}(\tilde u))^{-1} (I+\frac{\tau}{4} \mathcal{A}_{ir}(\tilde u)) u^{6/8} \\
u^{1} & = & (I-\frac{\tau}{4} \mathcal{A}_{ic}(\tilde u))^{-1} (I+\frac{\tau}{4} \mathcal{A}_{ic}(\tilde u)) u^{7/8} \\
\end{array}
\label{eq:marchuk-split}
\end{equation}
where $u^0 = u(t)$, and $u^1 = u(t+\tau)$.
\begin{theorem}[\cite{marchuk-split-1995}]\label{thm:marchuk-split}
Suppose each $\mathcal{A}_i(u) \le 0$ and sufficiently smooth w.r.t $u$, where $i \in \{ed,ec,ic,ir\}$.
Then the scheme \eqref{eq:marchuk-split} provides the second order of accuracy and is absolutely stable.
\end{theorem}

In each step of the splitting, a certain basic scheme is applied.
In fact, the Crank-Nicolson propagators in \eqref{eq:marchuk-split} can be substituted by any second-order accurate and stable scheme.
Furthermore, each physical process, e.g. ``electron convection'' may be further split componentwise.
Let us describe them one-by-one.

\subsubsection{Diffusion term for electron density}\label{sec:diff}
Due to the periodic boundary conditions, the diffusion matrix $\Delta$ obtained as a standard 5-point finite difference stencil for the Laplace operator in \eqref{eq:electron} is a two-level circulant,
which is easily diagonalizable by the two-dimensional Fourier transform.
Therefore, it is not difficult to apply any matrix function of $\Delta$ to a vector, for example, the exponential,
$$
\Delta = (F \otimes F) \mathrm{diag}(\lambda) (F^* \otimes F^*), \quad \exp(\tau \Delta)n_e = (F \otimes F) \mathrm{diag}(\exp(\tau\lambda)) (F^* \otimes F^*) n_e.
$$
Note that we need only four Fourier transforms of complexity $\mathcal{O}(n_x^2 \log n_x)$ and one multiplication by a diagonal matrix with $\mathcal{O}(n_x^2)$ complexity to compute the exact time step for the diffusion.
Therefore, the ``$ed$'' steps in \eqref{eq:marchuk-split} are performed as
$$
n_{e}^{3/8} = \exp(\tau/2 \cdot \Delta)n_{e}^{2/8}, \quad n_{e}^{6/8} = \exp(\tau/2 \cdot \Delta)n_{e}^{5/8}.
$$
The 5-point finite difference scheme provides $\Delta \le 0$ required in Theorem \ref{thm:marchuk-split}.

\subsubsection{Advection term for electron density}
For the electron advection parts(``$ec$''), we employ the Mac-Cormac space-time second-order accurate discretization.
Though such an advection alone may produce spurious oscillations with nonsmooth solutions,
in our case we gain some stabilization from the diffusion, since the Peclet number is small for the grid sizes used.
As soon as the Courant condition is satisfied, the Mac-Cormac scheme can be used inside the splitting \eqref{eq:marchuk-split}.

\subsubsection{Convection term for ion distribution}\label{sec:ic}
The convections in the ion equation \eqref{eq:vfp} (``$ic$'') possess one important property.
In all four terms, the convection velocity is independent on the variable, in which the gradient is taken.
In other words, each one-variate drift may be considered as a set of convections with constant velocities.
This problem can be solved analytically via the method of characteristics,
$$
f(x,y,v,w,t+\tau) = f(x-\tau v,y,v,w,t),
$$
and similarly for $y,v,w$.
For brevity, we present the procedure for $x$ only, repeating it for other variables in the splitting sense.
On a discrete level, the shift along a characteristic is substituted by a 5-point interpolation, the so-called ``cross'' scheme, exploited also in \cite{kovalev-plasma-2008,kovalev-plasma-2009}:
$$
f(x_i,t+\tau) = \alpha_{-2} f(x_{i-2},t) + \alpha_{-1} f(x_{i-1},t) + \alpha_0 f(x_i,t) + \alpha_1 f(x_{i+1},t) + \alpha_2 f(x_{i+2},t),
$$
where
\begin{equation}
\begin{array}{rclrcl}
\alpha_{-2} & = & \dfrac{c(c + 1)(c - 1)(c + 2)}{24}, & \alpha_{2} & = & \dfrac{c(c + 1)(c - 1)(c - 2)}{24}, \\
\alpha_{-1} & = & \dfrac{-c(c + 1)(c + 2)(c - 2)}{6}, & \alpha_{1} & = & \dfrac{-c(c - 1)(c + 2)(c - 2)}{6}, \\
\alpha_0 & = & \dfrac{(c + 1)(c - 1)(c + 2)(c - 2)}{4}, & c & = & \dfrac{\tau v}{h},
\end{array}
\label{eq:cross-formulae}
\end{equation}
where
$$
h=x_i-x_{i-1} = \dfrac{L}{n_x}
$$
is the mesh interval for the space.
If $x_i$ is a boundary point, the neighbors $x_{\pm 1},x_{\pm 2}$ are taken from the other end according to the periodicity.
Therefore, we obtain a 5-diagonal circulant matrix w.r.t. $x$, but note that it is parametrized by $v$ (see additional details in Section \ref{sec:final-scheme}).
This scheme is again second-order accurate in both space and time (provided the Courant condition $c<1$ is satisfied), and moreover, is more stable against spurious oscillations than the Mac-Cormac one,
which is important for the Vlasov-Fokker-Planck equation, containing only the convection terms in $x,y$ variables.
Note also that for $f^{0} \xrightarrow{\tau/2} f^{1/8}$ and $f^{7/8} \xrightarrow{\tau/2} f^{1}$ steps in \eqref{eq:marchuk-split} we need to take $c=\frac{\tau v}{2h}$ calculating \eqref{eq:cross-formulae}.

\subsubsection{Reaction term for ion distribution}\label{sec:ir}
Finally, let us see that the reaction step $\frac{\partial f}{\partial t} = f_0-f$ may be represented as an autonomous problem with approximately an orthogonal projector as a matrix, and hence admits again a fast but accurate time integration.
Note that the velocity part of $f$ is a perturbation to the Gaussian function, which has a band-limited Fourier spectrum.
It means that the simple rectangle quadrature formula for \eqref{eq:f0} becomes exact (up to the truncated part outside the domain) starting from some grid, when all harmonics in the spectrum are resolved.
So, $f_0$ at the grid points is computed as follows,
$$
\begin{array}{rcl}
f_0(x_{i_1}, y_{i_2}, v_{i_3}, w_{i_4}) & = & \left(\sum\limits_{j_3,j_4=1}^{n_v} f(x_{i_1}, y_{i_2}, v_{j_3}, w_{j_4}) h_v^2 \right) \\
& \cdot & \frac{1}{2\pi} \exp(-\frac{v_{i_3}^2}{2}) \exp(-\frac{w_{i_4}^2}{2}),
\end{array}
$$
where
$$
h_v = v_i - v_{i-1} = \dfrac{2 v_{max}}{n_v}
$$
is the mesh interval for the velocity.
This operation can be represented as the following matrix-by-vector product,
$$
f_0 = E f, \quad E = I_{n_x^2} \otimes \left(\frac{h_v}{\sqrt{2\pi}} \mathbf{e} \mathbf{1}^\top\right) \otimes \left(\frac{h_v}{\sqrt{2\pi}} \mathbf{e} \mathbf{1}^\top\right),
$$
where $\mathbf{e} = \begin{bmatrix}\exp(-\frac{v_i^2}{2})\end{bmatrix}_{i=1}^{n_v}$, and $\mathbf{1}$ is a vector of all ones.
But it is known that
$$
\frac{1}{2\pi} \int \exp(-\frac{v^2}{2}) \exp(-\frac{w^2}{2}) dv dw = 1,
$$
and henceforth,
$$
E^2=E + \O\left(\int_{|\mathbf{v}|>v_{max}} f(x,y,v,w) d\mathbf{v}\right), \quad f(x,y,v_{max},v_{max}) \sim \mathrm{e}^{-v_{max}^2}\ll 1.
$$
That is, $E$ is approximately an orthoprojector, as well as its complement $E_{\perp}=I-E$ in the reaction part
$$
\dfrac{\partial f}{\partial t} = (E-I) f = -E_{\perp} f,
$$
where $-E_{\perp} \le 0$.
Now, the exact solution of the latter problem reads $f(t+\tau) = \exp(-\tau E_{\perp}) f(t)$.
On the other hand, since for any degree $k$ it holds $E_{\perp}^k = E_{\perp}$, the exponential series simplifies to a scalar one,
$$
\exp(-\tau E_{\perp})  = I + \sum\limits_{k \ge 1} \dfrac{(-\tau)^k E_{\perp}^k}{k!} = I + E_{\perp} \sum\limits_{k \ge 1} \dfrac{(-\tau)^k}{k!} = I + E_{\perp} (\exp(-\tau)-1).
$$
That is, the propagation of the reaction part is performed exactly in time, and with the spectral accuracy in space as follows,
$$
f^{2/8} = \left(\mathrm{e}^{-\tau/2}I + (1-\mathrm{e}^{-\tau/2}) E\right) f^{1/8}, \quad  f^{7/8} = \left(\mathrm{e}^{-\tau/2}I + (1-\mathrm{e}^{-\tau/2}) E\right) f^{6/8}.
$$

\subsubsection{Poisson equation and reassembly of operators}
One equation not described yet is the computation of the electric potential \eqref{eq:poisson}.
However, in all steps of \eqref{eq:marchuk-split} except the first one, all matrices remain the same $\mathcal{A}(\tilde u)$.
On the other hand, $\phi$ is only needed in the construction of $\mathcal{A}(u)$.
That is, the Poisson equation needs to be solved only after the computation of $\tilde u$ in the very first step of the splitting process.
It can be done via the Fourier transform in the very similar way as described in Section \ref{sec:diff}.
After that, we reconstruct all matrices and proceed to the linear propagation steps $u^{1/8},\ldots,u^1$.

\subsubsection{Treating different timescales}
Due to different masses of electrons and ions, the timescales in Equations \eqref{eq:electron} and \eqref{eq:vfp} are also different.
In \cite{kovalev-plasma-2008,kovalev-plasma-2009} it was suggested to split one total time step (used for \eqref{eq:vfp}) into $20-40$ smaller subintervals in solution of the electron equation.
The main reason for that is the stronger Courant condition for the convection steps in \eqref{eq:electron} compared to those in \eqref{eq:vfp}.
In this work we can do the same: note that the splitting scheme \eqref{eq:marchuk-split} keeps the steps corresponding to the electron equations together ($u^{2/8} \rightarrow u^{6/8}$).
Therefore, we may further split them into $N_{ext} \sim 40$ steps without touching the ion part.
The same approach should be applied to the computation of $\tilde u$: we perform first the time step $\tau/2$ for the ion equation, and then $N_{ext}$ steps of size $\tau/(2 N_{ext})$ each for the electron part (the actual order of these steps does not matter, since only $\mathcal{O}(\tau)$ accuracy is required for $\tilde u$).

\subsection{Tensor formats and methods}\label{sec:tt}
The most difficult part of the problem described is the four-dimensional approximate Vlasov-Fokker-Planck equation for the ion distribution \eqref{eq:vfp}.
Given grid sizes $n_x$ in $x$ and $y$ directions, and $n_v$ in the velocity coordinates, the total amount of values to represent $f(t)$ scales as $n_x^2 n_v^2$, which rapidly goes beyond the available memory of a computer node.
To treat this problem, different approaches have been suggested.
In \cite{dimant-fb-pic-2004}, the particle method was used: we specify $f$ not at the uniform grid points, but at some sparsely (in fact randomly) distributed feasible amount of points, which can be moved along the domain.
Unfortunately, this Monte-Carlo-type setting provides a poor accuracy vs. number of points ratio: still a vast amount of data is needed to resolve the process with even qualitative correctness.
As an alternative, in \cite{kovalev-plasma-2008,kovalev-plasma-2009}, the deterministic ``brute force'' grid-based modeling on a high performance supercomputer was proposed.

In our work we also propose the uniform grid framework, but avoid the direct storage of all points.
Instead, we will approximate the discrete representations of functions by a smart sum of products of univariate items.

Given a function $f=f(x_1,\ldots,x_d)$, assume a $n$-point grid is introduced in each variable, $x_k \in \begin{bmatrix}x_k(i_k)\end{bmatrix}_{i_k=1}^n$.
Then the samples of $f$ at all grid points may be gathered to a $d$-dimensional array, or \emph{tensor} $f(i_1,\ldots,i_d)$.
Obviously, to store all the samples, one needs $n^d$ memory cells, which grows rapidly with both $n$ and $d$.
However, assume that $f$ can be written as a product of univariate functions, i.e.
\begin{equation}
f(x_1,\ldots,x_d) = f^{(1)}(x_1)  f^{(2)}(x_2) \cdots f^{(d)}(x_d).
\label{eq:rank1}
\end{equation}
Now, we may store only $n$ samples of each $f^{(k)}(x_k)$, since they define actually all values of $f$, obtaining in total $nd \ll n^d$ storage cost.

Of course, the ultimate separation \eqref{eq:rank1} occurs very rarely in practice, and an \emph{approximation} $f^{(1)}(x_1)  f^{(2)}(x_2) \cdots f^{(d)}(x_d) \approx f(x_1,\ldots,x_d)$ would generally lead to a poor accuracy.
To understand how can we generalize \eqref{eq:rank1}, consider first the two-dimensional case.
Since $F=\begin{bmatrix}f(i_1,i_2)\end{bmatrix}$ is a matrix, the variables are separated in a unique way, namely,
$$
\begin{array}{lll}
f(i_1,i_2) = \sum\limits_{\alpha=1}^n f^{(1)}_{\alpha}(i_1) f^{(2)}_{\alpha}(i_2) & \Leftrightarrow & F = F^{(1)} (F^{(2)})^\top.
\end{array}
$$
The approximation is introduced by shrinking the range of $\alpha$ from $n$ to some $r<n$ values.
The storage cost for $f^{(1)}$ and $f^{(2)}$ is now $2nr$, and we may ask about the dependence of the \emph{rank} $r$ on the desired accuracy threshold, $r=r(\epsilon)$.
What is important, the \emph{optimal} dependence is attained at the very particular decomposition, which is robustly computable using the well-established software like the LAPACK library.
\begin{theorem}
Any matrix $F$ admits the \emph{singular value decomposition} (SVD),
$$
f(i_1,i_2) = \sum\limits_{\alpha=1}^n U_{\alpha}(i_1) \sigma_{\alpha} V_{\alpha}(i_2),
$$
where $\sigma_1 \ge \sigma_2 \ge \cdots \ge \sigma_n \ge 0$, $\sigma_i^2 \in \lambda(F^*F)$, and $U$ and $V$ are matrices with orthonormal columns.
Moreover, the \emph{truncated} rank-$r$ SVD yields an optimal rank-$r$ approximation, i.e.
$$
\left\|F-\sum\limits_{\alpha=1}^r U_{\alpha} \sigma_{\alpha} V_{\alpha}^\top \right\|_2 = \min\limits_{\mathrm{rank}(G) = r} \|F-G\|_2 = \epsilon.
$$
\end{theorem}

Though generally $r(\epsilon)$ is unknown, if $f$ comes from the discretization of a smooth function, reasonable bounds may be proven, typically of the form \cite{tee-tensor-2003,tee-kron-2004}
$$
r = \O(\log^\beta (1/\epsilon) \log^\gamma(n)), \quad \beta,\gamma>0
$$

A generalization to the higher-dimensional case is not obvious.
The first idea that comes in mind is to take a simple sum of rank-1 components \eqref{eq:rank1} in the same way as is done in two dimensions.
This representation \emph{format} is known since \cite{hitchcock-sum-1927} under the names \emph{Canonical} Polyadic (CP) or PARAFAC.
It was then used heavily for the structuring of data, see e.g. the review \cite{kolda-review-2009} and references therein.
However, there is no SVD-like algorithm to compute the CP format, and moreover, it may suffer from an intrinsic instability \cite{desilva-2008}, which prevents efficient calculations.

Therefore, a better approach is to construct nested two-dimensional decompositions via e.g. the SVD algorithm.
For example, for $d=3$ we begin with the two-dimensional decomposition
$$
f(i_1, i_2, i_3) = \sum\limits_{\alpha_1} f^{(1)}_{\alpha_1} (i_1) g^{(1)}_{\alpha_1}(i_2, i_3).
$$
Proceeding in the same way, we may consider $g^{(1)}_{\alpha_1}(i_2, i_3)$ as a two-dimensional array
parameterized by $\alpha_1$.
However, this would significantly inflate the total number of the summation terms.
To avoid this, we should be keen to rewrite $g^{(1)}$ in the following way:
$$
g^{(1)}_{\alpha_1}(i_2, i_3) = F^{(2)}(\overline{\alpha_1 i_2}, i_3),
$$
where $\overline{\alpha_1 i_2}$ is a multi-index composed of $\alpha_1$ and $i_2$.
Then we may write the next two-dimensional decomposition in the form
$$
F^{(2)}(\overline{\alpha_1 i_2}, i_3) = \sum\limits_{\alpha_2} f^{(2)}_{\alpha_1 \alpha_2}(i_2) f^{(3)}_{\alpha_2}(i_3).
$$
Finally,
$$
f(i_1, i_2, i_3) = \sum\limits_{\alpha_1} \sum\limits_{\alpha_2}
f^{(1)}_{\alpha_1} (i_1) f^{(2)}_{\alpha_1 \alpha_2} (i_2) f^{(3)}_{\alpha_2}(i_3).
$$
Similarly, for an arbitrary $d$ we obtain
\begin{equation}
f(i_1,\ldots,i_d) = \sum\limits_{\alpha_1,\ldots,\alpha_{d-1}=1}^{r_1,\ldots,r_{d-1}}  f^{(1)}_{\alpha_1}(i_1) f^{(2)}_{\alpha_1,\alpha_2}(i_2) \cdots f^{(d-1)}_{\alpha_{d-2},\alpha_{d-1}}(i_{d-1}) f^{(d)}_{\alpha_{d-1}}(i_d).
\label{eq:tt}
\end{equation}
This format was called a \emph{Tensor Train}, or simply TT format \cite{ot-tt-2009,osel-tt-2011}, since each summation index $\alpha_k$ enters only the $k$ and $k+1$ blocks.
We see that we need to store only $d$ three-dimensional \emph{factors}, or \emph{TT blocks} $f^{(k)}$.
Provided the \emph{TT ranks} $r_k$ are all bounded by a moderate constant $r$, we obtain the memory estimate $\O(d n r^2)$.

One should note that such a structure appeared about 20 years ago in the community of quantum physics.
If we fix the indices $i_1,\ldots,i_d$ in \eqref{eq:tt}, the factors become \emph{matrices} with column resp. row indices $\alpha_{k-1},\alpha_k$, depending on $i_k$ as parameters.
Therefore, we may say that each element $f(i_1,\ldots,i_d)$ is equal to a product of $d$ matrices.
Since \cite{white-dmrg-1993}, this structure was used to handle high-dimensional wavefunctions of quantum states, and hence in \cite{klumper-mps-1993} it was named the \emph{Matrix Product States}.

A very important property of this structure is that it allows efficient counterparts of many linear algebra operations.
Linear combinations, scalar and pointwise products may be computed exactly with $\O(dnr^{p})$ complexity.
The matrix product is also available, but the structure of a \emph{multilevel} matrix requires some comments.
Recalling the definition of the Kronecker product $A = A^{(1)} \otimes A^{(2)}$, which means elementwise
$$
A(i_1,i_2,j_1,j_2) = A^{(1)}(i_1,j_1) A^{(2)}(i_2,j_2),
$$
we introduce the \emph{matrix} TT format, or the \emph{Matrix Product Operator} as the following generalization,
\begin{equation}
A(i_1,\ldots,i_d,j_1,\ldots,j_d) = \sum\limits_{\alpha_k=1}^{r_k} A^{(1)}_{\alpha_1}(i_1,j_1) A^{(2)}_{\alpha_1,\alpha_2}(i_2,j_2) \cdots A^{(d)}_{\alpha_{d-1}}(i_d,j_d).
\label{eq:ttm}
\end{equation}
Now, the MatVec operation involves a sequence of $d$ products of $A^{(k)}$ and $f^{(k)}$.

However, most of such operations return the result with larger TT ranks than that of the inputs.
In many cases they are not optimal, and a \emph{re-approximation} (or rounding) is required, i.e. reduction of the ranks to quasi-optimal ones for a given accuracy tolerance $\epsilon$.
We will denote this operation as follows,
\begin{equation}
\tilde x = \mathcal{T} x, \quad \mbox{or }\tilde x = \mathcal{T}_{\epsilon} x,
\label{eq:round}
\end{equation}
if we would like to stress the error threshold.
Fortunately, such a procedure can be performed using only $d-1$ QR and SVD decompositions and $2d-2$ matrix products of the total complexity $\O(dnr^3)$.
For more details we may refer to several books and surveys from the numerical mathematics \cite{khor-survey-2011,hackbusch-2012,larskres-survey-2013}, as well as the quantum physics \cite{PerezGarcia-mps-2007,schollwock-review-2011} communities.

We have mentioned only explicit operations, which can be performed in a finite number of steps with a guaranteed result.
Not of less importance are the \emph{implicit} solutions, which can be computed only iteratively and approximately in most cases (especially in such indirect representations as tensor product formats).
Though the tensor linear algebra equipped with the rounding procedure allows to think in terms of classical algorithms, use of the very special \emph{polylinear} structure of the TT product has appeared to be more efficient.

Given a functional $J(f)$ to minimize, we seek for $f$ in the TT format, and subsequently restrict the optimization to each TT block $f^{(k)}$, running through $k=1,\ldots,d$ during the iterations.
This approach was called the ALS (Alternating Least Squares, or later Alternating Linear Scheme).
Remarkable, if $J(f)$ is quadratic, the restriction $J_k(f^{(k)})$ remains quadratic, too, and we are left with a simple one-dimensional optimization in each step.
Unfortunately, such a straightforward realization does not allow to adapt the ranks at runtime, and moreover, is likely to converge slowly and stuck in spurious local minima of $J$.

As some remedy, the optimization over \emph{two} blocks $f^{(k)}f^{(k+1)}$ at a time was proposed firstly in the physics community under the name \emph{Density Matrix Renormalization Group} (DMRG) \cite{white-dmrg-1993} to compute the ground states (extreme eigenpairs) of spin systems, and then was brought to the numerical mathematics in \cite{holtz-ALS-DMRG-2012}.
After that, several improvements have been developed for the solution of linear systems \cite{DoOs-dmrg-solve-2011} and eigenvalue problems \cite{khos-dmrg-2010,dkos-eigb2013pre}, as well as the \emph{iterative rounding} in the MatVec operation \cite{Os-mvk2-2011}, and the adaptive interpolation of a tensor from smartly chosen samples \cite{so-dmrgi-2011proc,sav-qott-2013pre}.

Finally, a lying-in-between method with the complexity of the ALS and the rank adaptivity of the DMRG was proposed for linear systems \cite{ds-amr1-2013,ds-amr2-2013}.
This family was called AMEn (Adaptive Minimal Energy) due to its relation to the classical variational techniques, and the global convergence rate bound.

The mentioned iterative rounding operation is especially interesting for us, since we use only explicit time propagation schemes for the ion distribution $f$.
Given $f$ in \eqref{eq:tt} with the TT ranks bound $r(f)$, and $A$ in \eqref{eq:ttm} with the rank bound $r(A)$, the straightforward TT representation of the MatVec $g=Af$ is given with the rank bound $r(A) r(f)$,
while the optimal ranks are usually significantly smaller, $r(g) \sim r(f)$.
In principle, one could apply the TT rounding directly to the large-rank tensor $Af$, but the complexity overhead will be significant, especially if $r(A)$ is large, since
$$
\mathtt{work}_{direct}=\O(dn (r(A)r(f))^3) = \O(dn r(f)^3 \cdot r(A)^3).
$$

On the other hand, formulate an optimization problem $\min \|g-Af\|^2$, $g$ is sought in the form \eqref{eq:tt}.
The extremal condition for calculating a block $g^{(k)}$ can be satisfied with the cost of a product $g^\top A f$, which scales as
$$
\mathtt{work}_{g^\top A f} = d \cdot \O\left(n r(g) r(f)^2 r(A) + n^2 r(g) r(f) r(A)^2 + n r(g)^2 r(f) r(A) \right).
$$
The point is that $r(g)$ may be kept quasi-optimal during all iterations, for example, in a time propagation it is reasonable to initialize the algorithm with $f$.
Note that this approximation problem may be reformulated as a special case of a linear system, $\|g-Af\|^2 = \|Ig-(Af)\|^2$, and hence a counterpart of the AMEn algorithm with the corresponding theoretical considerations may be derived.
A more comprehensive study of the AMEn-rounding technique will be given in a separate work.
Here we just point out that in all cases of the Farley-Buneman simulations, the AMEn initialized with $f$ from the previous time step converges to the required threshold $\epsilon$ in two iterations, being more efficient than the direct TT rounding even for $r(A)=2$.
Thus, it will be the method of choice in the rest of the paper, and we will denote its action $\mathcal{T} \cdot$, in the same way as in \eqref{eq:round}.

\subsection{Getting things together}\label{sec:final-scheme}
Since the electron \eqref{eq:electron} and Poisson \eqref{eq:poisson} equations are two-dimensional, they require no further description.
We will apply the tensor approximation to the solution of the ion equation \eqref{eq:vfp}, and there are several details how to achieve a performance gain.

Up to this moment, we have not mentioned the initial state for the system, i.e. $n_e(0)$ and $f(0)$.
The initial velocity distribution is taken to be Maxwellian,
\begin{equation}
f(x,y,v,w,0) = n_i(x,y,0) \cdot \frac{1}{2\pi}\exp(-\frac{v^2}{2}) \exp(-\frac{w^2}{2}),
\label{eq:f_init}
\end{equation}
and the concentrations are taken randomly at each grid point,
$$
n_e(x_{i_1},y_{i_2},0) = n_0+M \cdot \mathrm{rand}, \quad n_i(x_{i_1},y_{i_2},0) = n_0+M \cdot \mathrm{rand}, \quad \mathrm{rand} \in \mathcal{N}(0,1),
$$
where $n_0$ is the average concentration, $\mathcal{N}(0,1)$ is a standard normal distribution, and $M$ is the magnitude constant, chosen in such a way that $E_{add} = \frac{1}{10}E_0$.
This choice is based on two considerations.
First, we would like to have all possible harmonics presented in the system, so the white noise is a good candidate for the initial state.
Second, we mimic the way how the Farley-Buneman process actually develops from random perturbations of the concentrations of particles.
Note that $M$ varies with the grid refinement (and also for different seeds of $\mathrm{rand}$!), but the experiments show that the behavior of the system remains the same as soon as $\frac{E_{add}}{E_0}$ is fixed.

Due to a random-looking structure w.r.t. $x,y$ coordinates, they cannot be well separated in a TT format.
Fortunately, the velocity distribution remains to be only a slight perturbation to the Maxwellian one during the whole process.
Therefore, a good separability of $x,y$ and $v,w$ coordinates is expected.
As a result, the TT structure \eqref{eq:tt} is used in the following way,
\begin{equation}
f(i_1,i_2,i_3,i_4) = f^{(1)}(i_1,i_2) f^{(3)}(i_3) f^{(4)}(i_4).
\label{eq:f_tt}
\end{equation}
Note that the initial state \eqref{eq:f_init} possesses the rank-1 TT representation.

Consider the convection steps.
Like in Section \ref{sec:ic}, we start from the part $\frac{\partial f}{\partial t}+v \frac{\partial f}{\partial x}=0$.
Each coefficient in the cross-formula \eqref{eq:cross-formulae} acts in connection with a certain periodic shift matrix.
Denote
\begin{equation}
S_1 = \begin{bmatrix}
       0 & 1 \\
       & 0 & 1 \\
       & & \ddots & \ddots \\
       & & & 0 & 1 \\
       1 & 0 & \cdots & & 0
      \end{bmatrix}, \quad
S_2 = \begin{bmatrix}
       0 & 0 & 1 \\
       & 0 & 0 & 1 \\
       & & \ddots & \ddots & \ddots \\
       & & & 0 & 0 & 1 \\
       1 & 0 & & \cdots  & 0 & 0 \\
       0 & 1 & 0 & \cdots & & 0
      \end{bmatrix},
\end{equation}
and $S_{-2} = S_2^{\top}$, $S_{-1}=S_1^{\top}$, $S_0=I$.
Then the one-dimensional convection propagation matrix reads
$$
M_c = \alpha_{-2} S_{-2} + \alpha_{-1} S_{-1} + \alpha_0 S_0 + \alpha_1 S_1 + \alpha_2 S_2.
$$
However, in our case each $\alpha_i$ is parametrized by the velocity.
That is, the action of $\alpha_i$ is given by a diagonal matrix,
$$
\Lambda_i = \left[I \otimes I\right] \otimes \mathrm{diag}(\alpha_i(v)) \otimes I,
$$
and for the total propagator we obtain the rank-5 representation:
$$
\begin{array}{rcl}
M_x & = & \left[S_{-2} \otimes I\right] \otimes \mathrm{diag}(\alpha_{-2}(v)) \otimes I + \left[S_{-1} \otimes I\right] \otimes \mathrm{diag}(\alpha_{-1}(v)) \otimes I \\
& + & \left[S_{0} \otimes I\right] \otimes \mathrm{diag}(\alpha_{0}(v)) \otimes I \\
& + & \left[S_{1} \otimes I\right] \otimes \mathrm{diag}(\alpha_{1}(v)) \otimes I + \left[S_{2} \otimes I\right] \otimes \mathrm{diag}(\alpha_{2}(v)) \otimes I,
\end{array}
$$
where the Kronecker products w.r.t. the $x$ and $y$ variables (in square brackets) are expanded, since the spatial dimensions are not separated, but the rest ones are kept implicitly in the TT structure.
In the same way we construct the other convection matrices:
$$
\begin{array}{rcl}
M_y & = & \left[I \otimes S_{-2}\right] \otimes I \otimes \mathrm{diag}(\alpha_{-2}(w))  + \left[I \otimes S_{-1}\right] \otimes I \otimes \mathrm{diag}(\alpha_{-1}(w))  \\
& + & \left[I \otimes S_{0}\right] \otimes I \otimes \mathrm{diag}(\alpha_{0}(w)) \\
& + & \left[I \otimes S_{1}\right] \otimes I \otimes \mathrm{diag}(\alpha_{1}(w)) + \left[I \otimes S_{2}\right] \otimes I \otimes \mathrm{diag}(\alpha_{2}(w)),
\end{array}
$$
$$
\begin{array}{rcl}
M_v & = & \mathrm{diag}\left(\alpha_{-2}\left(V_v\right)\right) \otimes S_{-2} \otimes I   + \mathrm{diag}\left(\alpha_{-1}\left(V_v\right)\right) \otimes S_{-1} \otimes I  \\
& + & \mathrm{diag}\left(\alpha_{0}\left(V_v\right)\right) \otimes S_{0} \otimes I \\
& + & \mathrm{diag}\left(\alpha_{1}\left(V_v\right)\right) \otimes S_{1} \otimes I + \mathrm{diag}\left(\alpha_{2}\left(V_v\right)\right) \otimes S_{2} \otimes I,
\end{array}
$$
where $V_v = -\frac{\partial \phi}{\partial x}$, and
$$
\begin{array}{rcl}
M_w & = & \mathrm{diag}\left(\alpha_{-2}\left(V_w\right)\right) \otimes I\otimes S_{-2}
 +  \mathrm{diag}\left(\alpha_{-1}\left(V_w\right)\right) \otimes I\otimes S_{-1}   \\
& + & \mathrm{diag}\left(\alpha_{0}\left(V_w\right)\right) \otimes I\otimes S_{0}  \\
& + & \mathrm{diag}\left(\alpha_{1}\left(V_w\right)\right) \otimes I\otimes S_{1}
 +  \mathrm{diag}\left(\alpha_{2}\left(V_w\right)\right) \otimes I\otimes S_{2},
\end{array}
$$
where $V_w = \frac{e E_0}{m_i v_{T_i} \nu_{in}} - \frac{\partial \phi}{\partial y}$.
Now, the ``ion convection'' splitting steps (e.g. $u^{1/8}$) in \eqref{eq:marchuk-split} are performed as follows,
$$
f^{1/8} = \mathcal{T}  M_w \mathcal{T} M_v \mathcal{T}  M_y \mathcal{T} M_x f^0,
$$
where $\mathcal{T}$ denotes the (AMEn-)rounding operation in the TT format (see Section \ref{sec:tt}).

The reaction part is easier, since it involves only the rank-2 TT structure of the propagation matrix.
Indeed, as was shown in Section \eqref{sec:ir}, the reaction step may be computed as follows,
$$
f^{2/8} = \mathcal{T} M_r f^{1/8},
$$
where
$$
M_r = \mathrm{e}^{-\tau/2} I \otimes I \otimes I \otimes I + (1-\mathrm{e}^{-\tau/2}) I \otimes I \otimes \left(\frac{h_v}{\sqrt{2\pi}} \mathbf{e} \mathbf{1}^\top\right) \otimes \left(\frac{h_v}{\sqrt{2\pi}} \mathbf{e} \mathbf{1}^\top\right),
$$
and similarly for $f^{7/8}$.

\section{Numerical experiments}\label{sec:num}
We implement the prototype solution software in MATLAB, and run at the Intel CPU  @2GHz machine with 70Gb of shared memory.
The version with the TT representation of the ion distribution is compared with the one where all elements of $f$ are stored directly as a full array (\emph{full} format).
Since MATLAB is not very suited for parallelization, and in order to fit the full format simulation into memory, we had to conduct the experiments not with the real electron mass, but with the modified parameters, chosen according to \cite{kovalev-plasma4effect-2009,kovalev-plasma-2009}, see Table \ref{tab:param_mmass} (note that the parameters are given already in the rescaled quantities according to Table \ref{tab:rescale}).
\begin{table}[h!]
\centering
\caption{Simulation parameters corresponding to the modified electron mass}
\label{tab:param_mmass}
\begin{tabular}{|c|cc|c|c|cc|} \hline
\multicolumn{3}{|c|}{Physical parameters} & \quad \quad & \multicolumn{3}{|c|}{Discretization parameters} \\ \hline
$T_i$ & 300$\cdot$1.3806505e-23 & [J]       &       & $L$ & $50$ &  \\ \hline
$T_e$ & 1 &                                 &       & $v_{max}$ & $6$ & \\ \hline
$E_0$ & 0.05 & [V/m]                        &       & $n_x$ & $250$ & \\ \hline
$B_0$ & 5e-5  & [T]                         &       & $n_v$ & $31$ & \\ \hline
$n_0$ & 1e+10  & [$m^{-3}$]                 &       & $\tau$ & $0.01$ & \\ \hline
$m_i$ & 4.9936722e-26 & [kg]                &       & $N_{ext}$ & $40$ & \\ \hline
$\nu_{in}$ & 1800 & [1/s]                   &       & $\epsilon$ & $0.05 \cdot \frac{|f(t)-f(t-\tau)|}{|f(t-\tau)|}$ & \\ \hline
\multicolumn{3}{|c|}{Model scales}          &       & \multicolumn{3}{|c|}{Physical constants} \\ \hline
$\psi$ & $0.1575$ &                         &       & $e$ & 1.60217653e-19 & [C] \\ \hline
$\theta$ & $0.03528$ &                      &       & $\eps_0$ & 8.85418781e-12 & [$\frac{F}{m}$] \\ \hline
$l$ & $0.16$ & [m]                          &       & $m_e$ & 3.97950489e-29 & [kg] \\ \hline
\end{tabular}
\end{table}
This allows us to use smaller grid sizes (due to a smaller Peclet number) and accelerate the simulation without affecting the qualitative behavior of the system, hence focusing on the demonstration of the tensor structuring potential rather than difficulties arising from the turbulent nature of the electron equation \eqref{eq:electron}.
Several approaches to treat this drawback are discussed in the conclusion and planned for a future work.

\begin{figure}[h!]
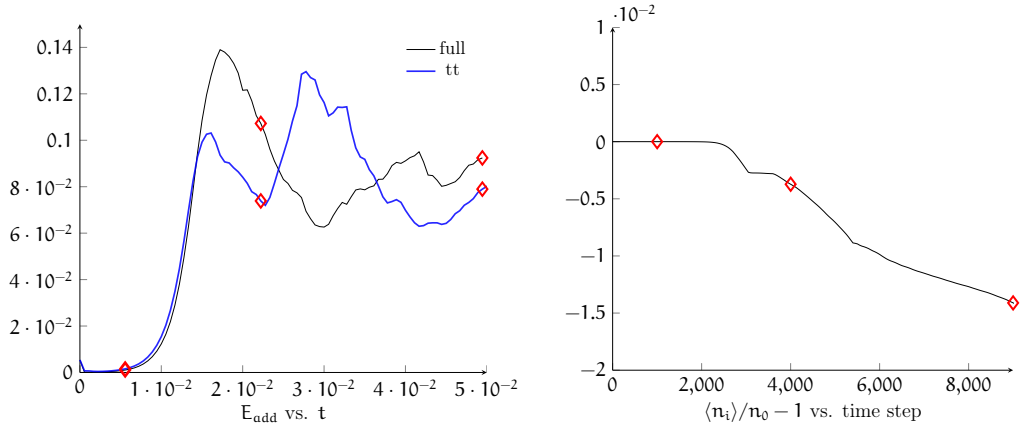

\centering
\resizebox{.52\textwidth}{!}{\input{./Pic/pgfart.sty} \input{./Pic/Eext.tikz}} \hfil
\resizebox{.46\textwidth}{!}{\input{./Pic/pgfart.sty} \input{./Pic/errmean.tikz}}

\caption{Additional electric field in the full format and TT-based simulations (left). Error of the average ion concentration in the TT method (right).}
\label{fig:Eext}
\end{figure}

\begin{figure}[h!]
\centering
\includegraphics[width=0.49\linewidth]{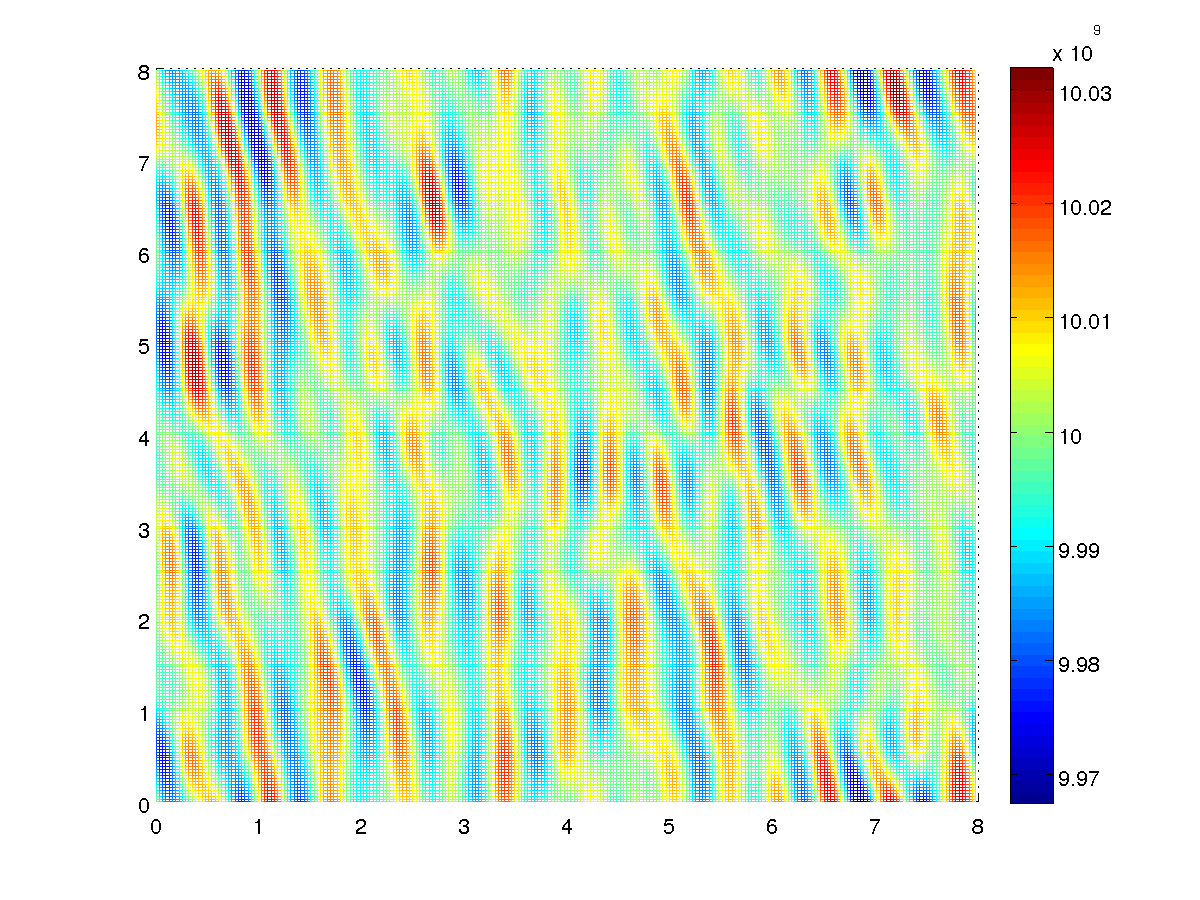} \hfil
\includegraphics[width=0.49\linewidth]{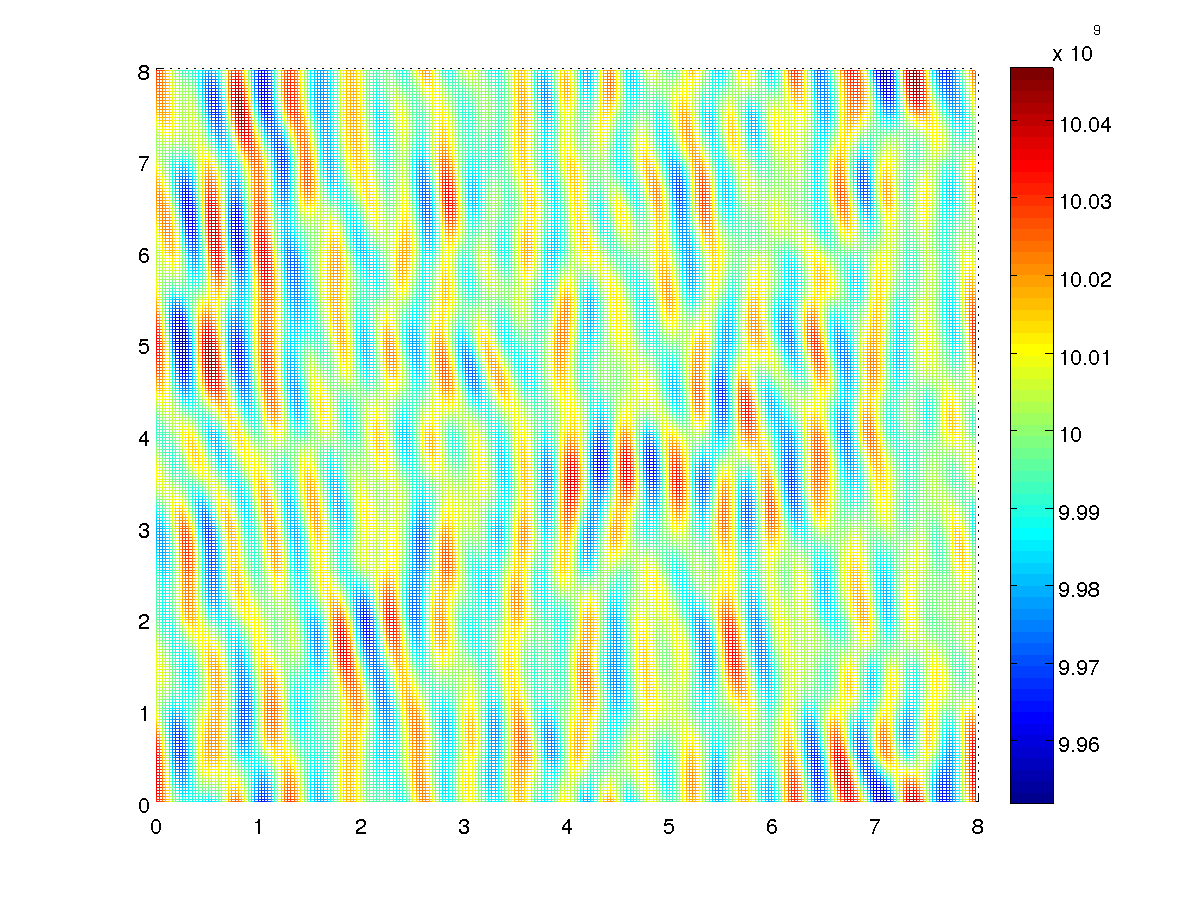} \\
\includegraphics[width=0.49\linewidth]{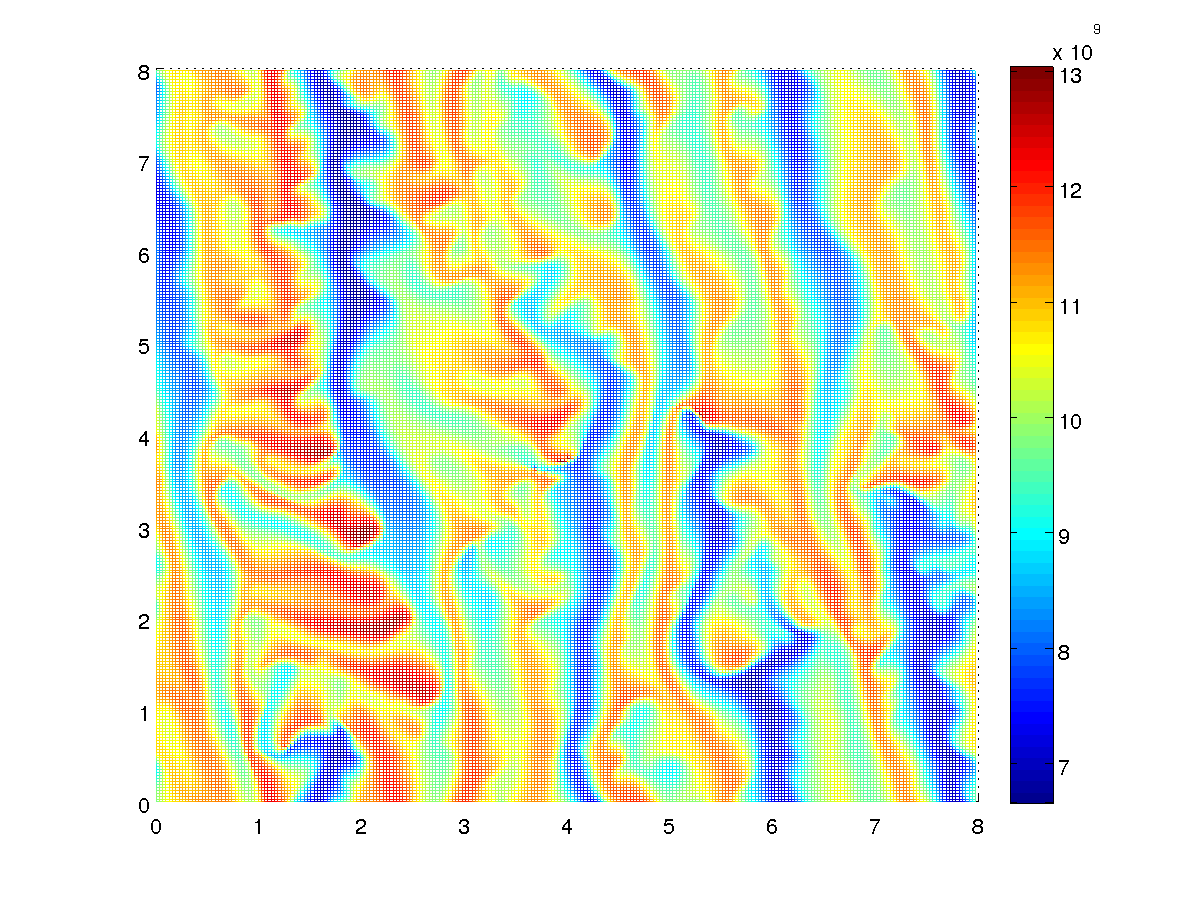} \hfil
\includegraphics[width=0.49\linewidth]{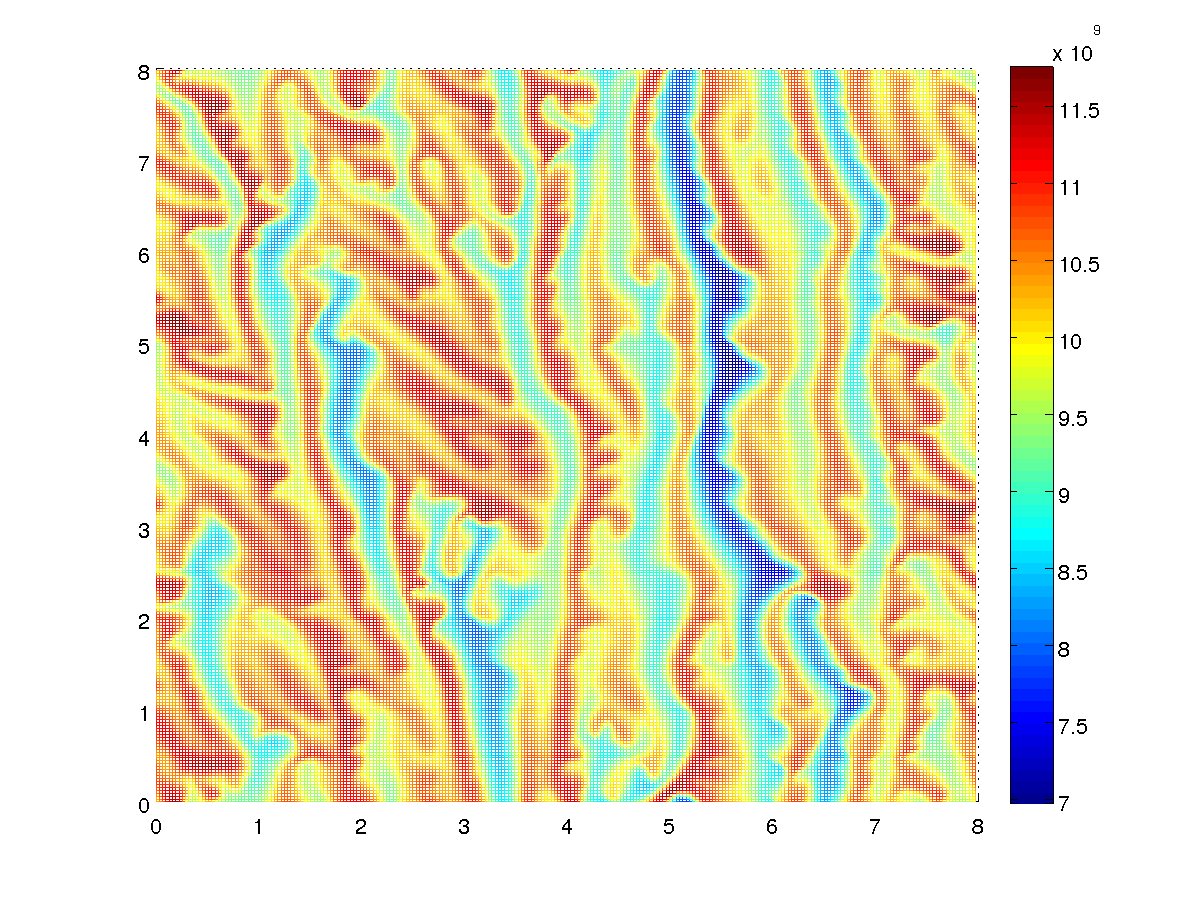} \\
\includegraphics[width=0.49\linewidth]{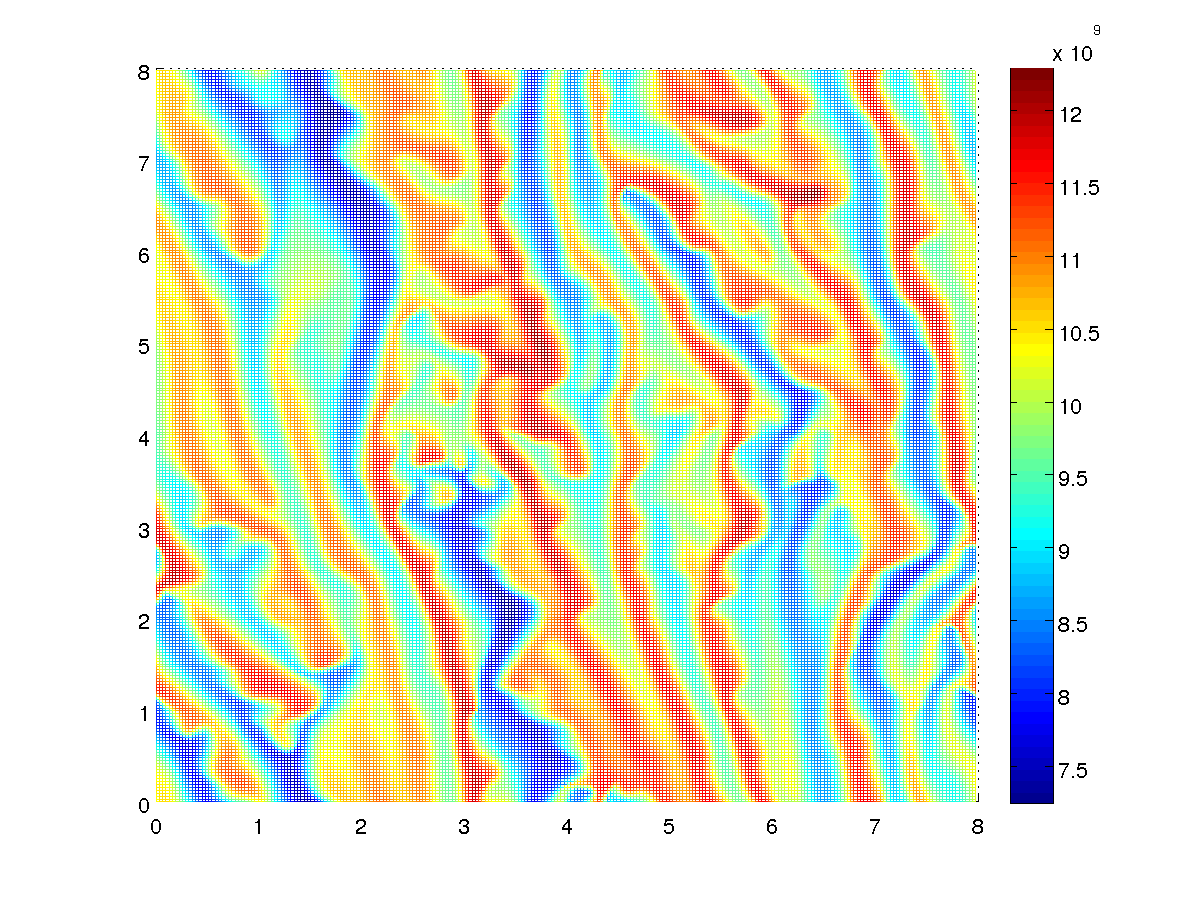} \hfil
\includegraphics[width=0.49\linewidth]{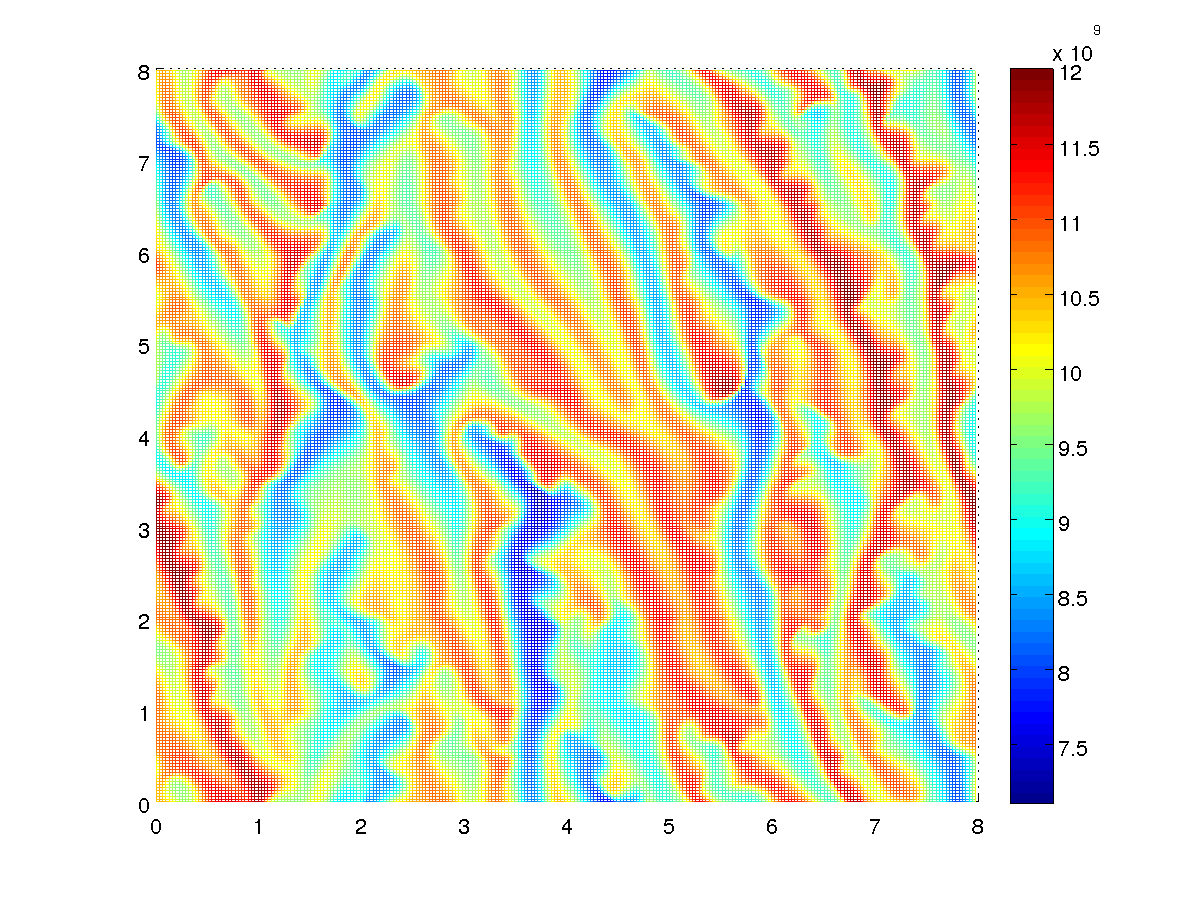} \\

\caption{Electron concentrations in the full format (left) and TT-based (right) simulations at different time steps: 1000 (top), 4000 (middle) and 9000 (bottom). Axes: $x$ and $y$ [m].}
\label{fig:ne}
\end{figure}

\begin{figure}[h!]
\centering
\includegraphics[width=0.49\linewidth]{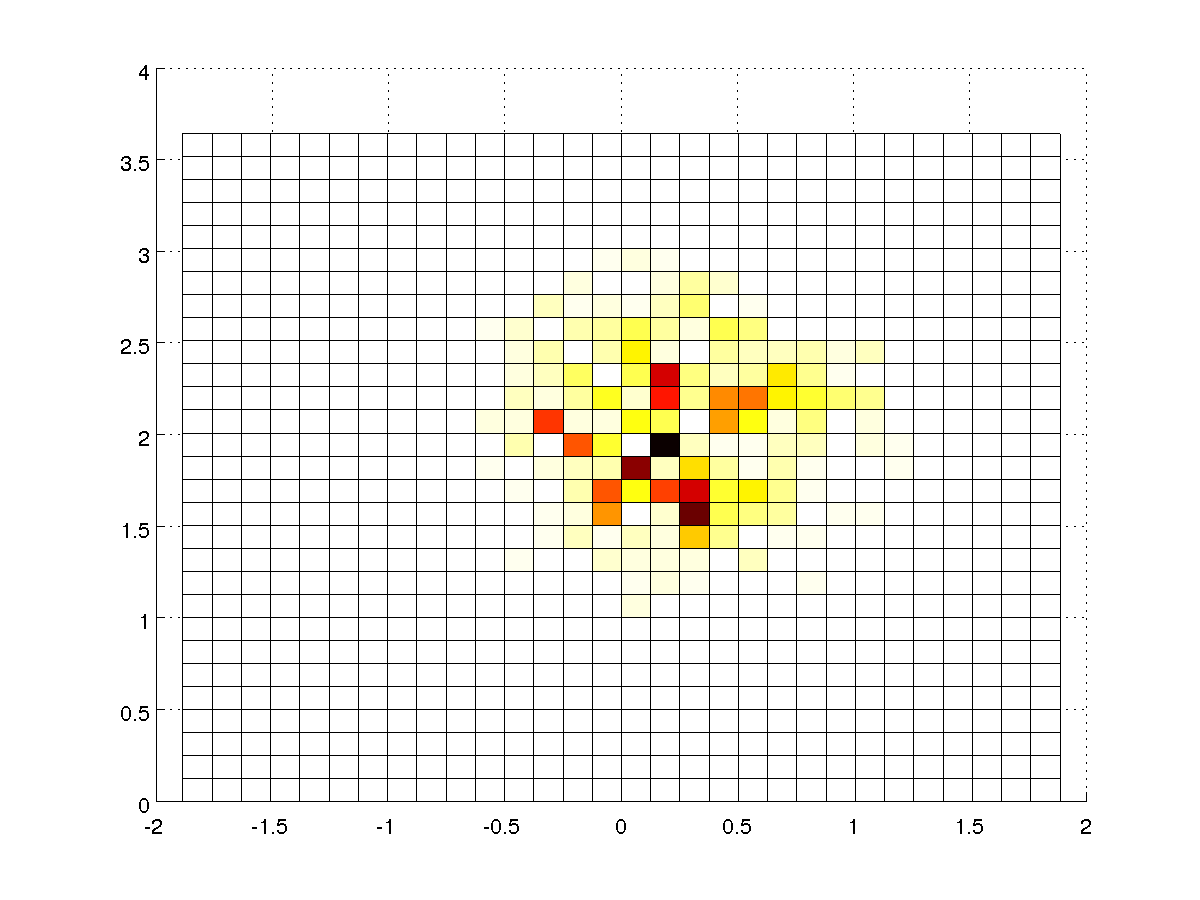} \hfil
\includegraphics[width=0.49\linewidth]{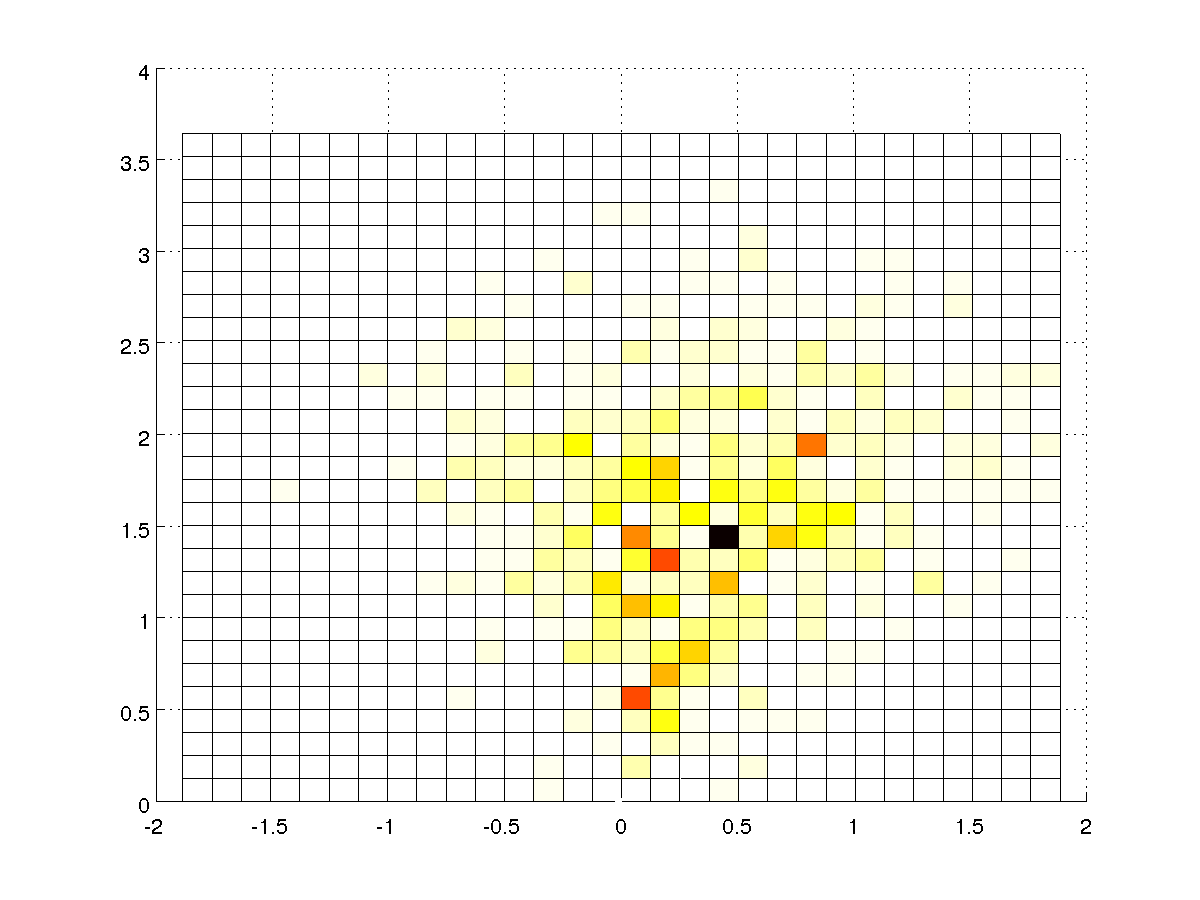}

\caption{Spectral intensities $\hat E^2$ of the electric field in the TT-based simulations at time steps: 1000 (left), and 9000 (right).  Axes: $k_x$ and $k_y$ [1/m].}
\label{fig:fE2}
\end{figure}

\begin{figure}[h!]
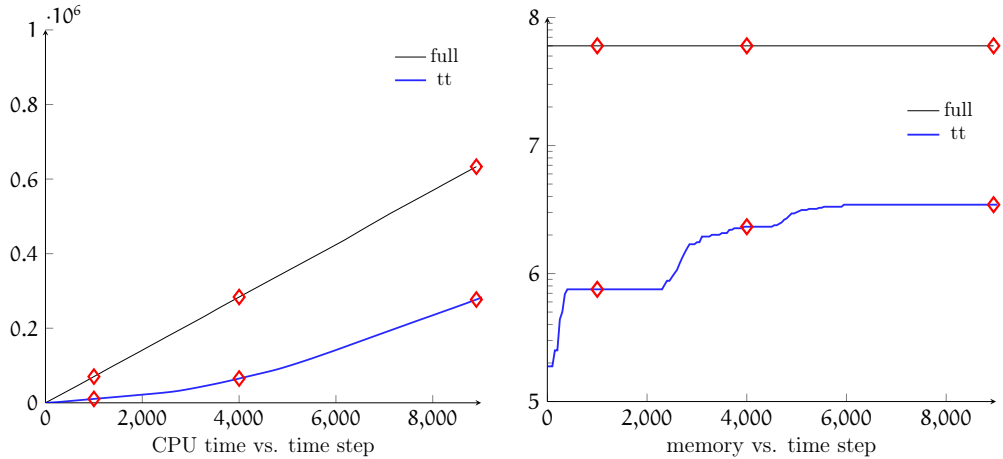

\centering
\resizebox{.48\textwidth}{!}{\input{./Pic/pgfart.sty} \input{./Pic/ttimes.tikz}} \hfil
\resizebox{.48\textwidth}{!}{\input{./Pic/pgfart.sty} \input{./Pic/mem.tikz}}
\caption{Cumulative CPU time, sec. (left) and the storage of $f$ (right) in the full format and TT-based simulations.}
\label{fig:ttimes}
\end{figure}

So, we perform 9000 time steps, corresponding to $t = 0.05$ s.
Due to a qualitative nature of the model (BGK collision operator, modified electron mass), we may admit rather high tensor rounding error $\epsilon$, keeping the correction at each time step with $5\%$ accuracy (see Table \ref{tab:param_mmass}).
Therefore, an accuracy of about $10\%$ may be expected for the output quantities.
Considering the average additional field $E_{add}$ (Figure \ref{fig:Eext}, left), and the average concentration of ions (Figure \ref{fig:Eext}, right), we conclude that this accuracy level is maintained.
Note that since the tensor rounding performs in fact orthogonal projections of the solution, the norm of the solution decreases from step to step, which can be clearly seen in the degeneracy of the average concentration $\langle n_i \rangle$.

Looking at the electric field $E_{add}$ (Figure \ref{fig:Eext}), as well as the electron concentrations (Figure \ref{fig:ne}) at different time steps in more details, we may see that during the initial stage of the development of the Farley-Buneman process, the full format and approximate TT solutions coincide with high precision (time steps $\le 2000$).
The same holds true for the timescales of the nonlinear saturation (Table \ref{tab:errs_E}, last row): the system is considered to enter a substantially nonlinear stage if the additional electric field begins to decrease after the linear growth (note that there is also a small region of oscillating field in the very beginning of the process).

However, a nonlinear system becomes more sensitive to the perturbations arising from the tensor truncations,
and the solutions develop in significantly different ways during the transient processes (time step $\sim 4000$).

Finally, the fully saturated nonlinear system enters a stationary region (time steps up to $9000$), where the additional field oscillates around its average value.
Note that despite significantly different distributions of concentrations (Figure \ref{fig:ne}), the statistical quantities are much less sensitive to the solution errors, see Table \ref{tab:errs_E}.
Considering the spatial harmonics of the electric field (Figure \ref{fig:fE2}), we observe that in the beginning of the process, the spectrum is almost isotropic w.r.t. the $x$-axis, while after the saturation anisotropic components appear, i.e. the model predicts the rotation of the drift vector correctly.

\begin{table}[h!]
\centering
\caption{Statistical electric field outputs}
\label{tab:errs_E}
\begin{tabular}{|c|c|c|c|} \hline
& Full solver & TT solver & Error \\ \hline
$\frac{1}{0.05-0.03}\int\limits_{0.03}^{0.05}E_{add}(t) dt$ & 8.3083e-02 & 8.0523e-02 & $3.08\%$ \\ \hline
$\max\limits_{0.015 \le t \le 0.03} E_{add}(t)$ & 1.3903e-01 & 1.3200e-01 & $5.05\%$ \\ \hline
$\min\limits_{t > 0.01} t:~\frac{d E_{add}(t)}{dt}<0$ & 1.7267e-02 & 1.5989e-02 & $7.40\%$ \\ \hline
\end{tabular}
\end{table}

Being assured with the correctness of the solution, consider the computational complexity, Figure \ref{fig:ttimes}.
The TT ranks grow during the development of the Farley-Buneman process, and stabilize at the maximal value $55$ after the saturation of the system (instead of ranks, we present directly the amount of memory cells to store each $f(t)$ in Figure \ref{fig:ttimes}, right).
Most storage and computational cost is due to the first TT block $f^{(1)}(i_1,i_2)$ in \eqref{eq:f_tt}, which contains $n_x^2 r$ elements.
So, $r=55$ should be compared with $n_v^2=961$ in the full format, thus giving the memory reduction factor greater than $17$.

The computational complexity grows faster with ranks than the storage one.
Therefore, there is not so impressive acceleration as for the memory consumption.
Nevertheless, the solution time of the TT solver is at least $2$ times smaller than of the full format one during the present simulation (Figure \ref{fig:ttimes}, left), and moreover, grows milder with the time steps, i.e. the difference would be even more significant if further simulation is needed.

The same prediction may be formulated with respect to the grid refinement.
In our case, we just could not fit all the data required for the full format solution into the memory when $n_x \gtrsim 300,~n_v \gtrsim 30$.
However, since the TT ranks are in most cases almost stable w.r.t. the grid sizes, the separation of variables is expected to be more efficient for finer grids.

\section{Conclusion}\label{sec:conclusion}
We have considered the tensor product approach to the adaptive black-box model reduction in the simulation of the two-component plasma driven by the electric fields.
As a particular example, we follow the hybrid kinetic-liquid model of the Farley-Buneman instability in two dimensions (perpendicular to the geomagnetic field), proposed in the previous papers \cite{dimant-fb-pic-2004,kovalev-plasma-2008,kovalev-plasma-2009}.
The most difficult part of the solution process is the calculation of the four-dimensional distribution function for ions, and it is the step where the approximate separation of variables (namely, the Tensor Train technique) is applied.

Since the electric forces introduce only a slight perturbation to the Maxwellian velocity distribution,
the spatial and velocity variables appear to be well separable.
Though some cost reduction possibility due to moderate TT ranks was already confirmed by our preliminary experiments,
an optimality of the actual rank values is doubtful.
Since we exploited only low-order discretization schemes, the solution (especially in presence of convections) is likely to be ``poisoned'' by a poorly structured discretization noise.
Depending on the actual magnitude, a noise can increase separations ranks far above the optimal ones of the exact function.
By employing high-order schemes we may expect a substantial increase in the accuracy, but with even smaller amount of memory as a bi-product.
The computational complexity of e.g. spectral element discretizations seems to be reasonable as well, since periodic boundary conditions yield circulant discrete counterparts of operators, which can be handled efficiently via the FFT.

Another way to gain speedups can be seen in the parallelization of the algorithms.
Contrarily to the previous work \cite{kovalev-plasma-2008}, we consider here only prototype MATLAB implementations.
However, since the spatial dimensions are still large even in the TT structure \eqref{eq:f_tt},
parallel versions of the TT algorithms run at high-performance machines may appear to be very efficient.

Finally, a question that may be addressed to physicists is how to improve the quality of the mathematical model itself.
For example, as soon as robust methods for handling the distribution functions in the whole phase space will be available, one may think of getting rid of simplified hybrid models, formulating the full Vlasov-Fokker-Planck equation for both electrons and ions.
All these questions are planned to be considered in a future research.

\bibliography{./bibtex/algebra,./bibtex/misc,./bibtex/our,./bibtex/tensor,./bibtex/dmrg}
\end{document}

%% file: Pic/Eext.tikz

\begin{tikzpicture}

 \begin{axis}[%
  xmode=normal,ymode=normal,
  cycle list name=fbm,
  xmin=0, xmax=0.05,
  ymin=0, ymax=0.15,
  yminorticks=true,
  legend style={at={(0.97,0.97)},anchor=north east}]

  \pgfplotsset{xlabel/.add={$E_{add}$ vs. $t$}}
  \addplot+[] table[header=false, x index=1, y index=2]{./dat/fb_Eext_full_250.dat}; \addlegendentry{full};
  \addplot+[] table[header=false, x index=1, y index=2]{./dat/fb_Eext_amen_250.dat}; \addlegendentry{tt};

  \addplot+[] table[header=false, x index=2, y index=3]{./dat/fb_stop_Eext_250.dat};
  \addplot+[] table[header=false, x index=2, y index=4]{./dat/fb_stop_Eext_250.dat};

\end{axis}
\end{tikzpicture}%

%% file: Pic/errmean.tikz

\begin{tikzpicture}

 \begin{axis}[%
  xmode=normal,ymode=normal,
  cycle list name=fbm,
  xmin=0, xmax=9000,
  ymin=-0.02, ymax=0.01,
  yminorticks=true,
  legend style={at={(0.97,0.97)},anchor=north east}]

  \pgfplotsset{xlabel/.add={$\langle n_i \rangle/n_0-1$ vs. time step}}
  \addplot+[] table[header=false, x index=0, y index=2]{./dat/fb_errmean_amen_250.dat};

  \pgfplotsset{cycle list shift=1}
  \addplot+[] table[header=false, x index=1, y index=3]{./dat/fb_stop_errmean_amen_250.dat};

\end{axis}
\end{tikzpicture}%

%% file: Pic/ttimes.tikz

\begin{tikzpicture}

 \begin{axis}[%
  xmode=normal,ymode=normal,
  cycle list name=fbm,
  xmin=0, xmax=9000,
  ymin=0, ymax=1e6,
  yminorticks=true,
  legend style={at={(0.97,0.97)},anchor=north east}]

  \pgfplotsset{xlabel/.add={CPU time vs. time step}}
  \addplot+[] table[header=false, x index=0, y index=2]{./dat/fb_ttimes_full_250.dat}; \addlegendentry{full};
  \addplot+[] table[header=false, x index=0, y index=2]{./dat/fb_ttimes_amen_250.dat}; \addlegendentry{tt};

  \addplot+[] table[header=false, x index=1, y index=3]{./dat/fb_stop_ttimes_250.dat};
  \addplot+[] table[header=false, x index=1, y index=4]{./dat/fb_stop_ttimes_250.dat};

\end{axis}
\end{tikzpicture}%

%% file: Pic/mem.tikz

\begin{tikzpicture}

 \begin{axis}[%
  xmode=normal,ymode=log,
  cycle list name=fbm,
  xmin=0, xmax=9000,
  ymin=1e5, ymax=1e8,
  yminorticks=true,
  legend style={at={(0.97,0.8)},anchor=north east}]

  \pgfplotsset{xlabel/.add={memory vs. time step}}
  \addplot+[] table[header=false, x index=0, y index=2]{./dat/fb_mem_full_250.dat}; \addlegendentry{full};
  \addplot+[] table[header=false, x index=0, y index=2]{./dat/fb_mem_amen_250.dat}; \addlegendentry{tt};

  \addplot+[] table[header=false, x index=1, y index=3]{./dat/fb_stop_mem_250.dat};
  \addplot+[] table[header=false, x index=1, y index=4]{./dat/fb_stop_mem_250.dat};

\end{axis}
\end{tikzpicture}%